\documentclass[12pt,a4paper]{article}
\usepackage[english,russian]{babel}
\usepackage[utf8]{inputenc}
\usepackage{amsmath}
\usepackage{amssymb}
\usepackage{amsfonts}
\usepackage{amsthm}
\usepackage{eucal}
\usepackage{amstext,upref}
\usepackage{geometry}
\usepackage{latexsym}
\usepackage{hyperref}%гиперссылки
\usepackage{xcolor}
\definecolor{linkcolor}{HTML}{799B03} % цвет ссылок
\definecolor{urlcolor}{HTML}{799B03} % цвет гиперссылок
\hypersetup{pdfstartview=FitH,  linkcolor=linkcolor,urlcolor=urlcolor, colorlinks=true}

\newcommand{\spb}{\par\vspace{2mm}\noindent}

\newcommand{\vbpb}{\par\vspace{1cm}\noindent}

\begin{document}
\begin{center}
	%title
{ \Large \bf Multivariate analog of the Le Roy-Lindel\"of theorem about
 power series analytic continuation} 
\end{center}
%name
\centerline{\large Aleksandr Mkrtchyan \footnote{This work was supported by the Science Committee of RA (Research project № PHDPD-23/1-04) and  supported by the Ministry of Science and Higher Education of the Russian Federation (agreement no. 075–15–2022–287).}}
\spb
%e-mail: diachenko@sfedu.ru
%affil
\centerline{Institute of Mathematics  NAS RA, Siberian Federal University, Saint Petersburg University   }
\spb

%text

%\centerline{\large A.\,J.\,Mkrtchyan 
%\footnote{This work was supported by the Science Committee of RA (Research project № PHDPD-23/1-04) and  supported by the Ministry of Science and Higher Education of the Russian Federation (agreement no. 075–15–2022–287).}}\\

{\bf Abstract.} We consider the problem of continuability into a sectorial domain for multiple power series and prove the multivariate analog of the Le Roy-Lindel\"of theorem, i.e. establish a connection between the growth of the  holomorphic function interpolating the series coefficients in the imaginary subspace and the  sectoral domain where the multiple power series is analytically continued to. \\


One of the methods of studying the question of analytical continuation of power series is interpolation of the coefficients of the series. With this approach 
Le Roy and Lindelöf \cite{Roy},\cite{L1} obtained conditions under which a series  analytically extends into a sector.  Note that the theorem, gave a connection between the sector and the growth of the interpolation function. More precisely, the type of interpolation function  must be less than $\pi$ on the closed half-plane $Re z \geq 0$.\\

\textbf{Theorem } (Le Roy, Lindel\"of) 
\textit{For a function } $\varphi$ \textit{holomorphic } \textit{in} $\mathbb{R}_{\geq0}+i\mathbb{R}$  \textit{and of exponential type} $\sigma < \pi$  \textit{the series} 
\begin{equation*} 
{f(z)=\sum_{n=0}^\infty \varphi(n) {z}^{n},     }   
\end{equation*}
\textit{admits an analytic continuation to the sector} $ \mathbb{C}\setminus \Delta_\sigma$. 
\textit{Here} $\Delta_\sigma  $ \textit{is the sector}  $\{z=re^{i\theta} \in \mathbb{C} : |\theta|\leq \sigma\}, \; \sigma \in [0,\pi) .$ \\

It should also be noted that for power series with unit radius of convergence Arakelyan obtained a criterion for the continuation  into a given sector, in which it is necessary to know the growth of the interpolation function in the positive open half-plane.\\

\textbf{Theorem } (Arakelian\cite{Ar84}) \textit{The sum of the series} 
\begin{equation*} 
{f(z)=\sum_{n=0}^\infty f_{n} {z}^{n}, \;\; \;  \limsup_{n \to \infty}\sqrt[n]{|f_{n}|}=1     }   
\end{equation*}
 \textit{ extends analytically to the open sector  }$ \mathbb{C}\setminus \Delta_\sigma$ \textit {if and only if there is an entire function}  $\varphi(\zeta)$    \textit{of exponential type interpolating the coefficients} $f_{n}$ \textit{ whose indicator function  } $h_\varphi(\theta)$ \textit {satisfies the condition}
\begin{equation*}{\label{1.33}}
h_{\varphi}(\theta)\leq \sigma |\sin \theta| \;\;  \textit{for} \;\; |\theta|< \frac{\pi}{2}.
\end{equation*}

Here $h_\varphi(\theta)$ is the indicator of an entire function of exponential type  introduced by Phragm\'en and Lindel\"of \cite{Bi}:
\begin{equation}\label{indicator}
h_{\varphi}(\theta):=\limsup_{r \to \infty}  \frac{\ln\left|f\left(re^{i\theta}\right)\right|}{r}, \;\; \theta \in \mathbb{R}.
\end{equation}

For multiple power series within this approach there are few results. In \cite{Mk3}, a criterion for continuability of a multiple power series across a boundary set of polyarcs was obtained in terms of asymptotic behaviour of  an entire function interpolating the coefficients of the series.
Also obtained conditions under which the multiple power series is analytically extends in to the sectoral domain.

In \cite{Mk2} we give a condition for continuability of a multiple power series into a sectorial domain defined by a piece-wise affine majorant for a function interpolating the series coefficients.
In \cite{Mk1} the condition of continuability into a sectorial domain  is given in terms of entire function interpolating the coeffcients of power
series.

Consider a multiple power series
\begin{equation}\label{1.0}
{f(z)=\sum_{k\in\mathbb N^n} f_{k} {z}^{k},   }
\end{equation}
where $z^k=z_1^{k_1}...z_n^{k_n}.$	
We say that a function $\varphi(\zeta)$ of $n$ complex variables  $\zeta=(\zeta_1,...,\zeta_n)$ \textit{interpolates} the coefficients of the series (\ref{1.0}) if
\begin{equation*}{\label{1.2}}
{ \varphi(k)=f_{k}\;\;\text{for all} \;\; k \in \mathbb N^n.}
\end{equation*}
The complex variables $\zeta_j$ we write as  $\zeta_j=\xi_j+i\eta_j$, thus $\xi \; $ is a vector of the real subspace of  $\mathbb{C}^n$, and  $\eta$  is a vector of the imaginary subspace. 

Introduce a piece-wise affine function
\begin{equation}\label{t1.0}{    \tilde{g}(\eta)= \sum_{p=1}^{q}\varepsilon_p|a_p^1\eta_1+...+a_p^n\eta_n+a_p^0| }
\end{equation}
with coefficients $ \varepsilon_p=\pm 1.$ Together with $\tilde g(\eta)$ we consider the following piece-wise linear function 
\begin{equation}\label{t1.01}{    {g}(\eta)= \sum_{p=1}^{q}\varepsilon_p|a_p^1\eta_1+...+a_p^n\eta_n|+\pi\sum_{k=1}^{n}(\eta_k-|\eta_k|).  }
\end{equation}
The last function has non-smooth points (points where linearity breaks) on the set of hyperplanes 
$$a_p^1\eta_1+...+a_p^n\eta_n=0, \,\, p=1,...,q \;\; \text{and} \;\; \eta_k=0, \,\, k=1,...,n.$$
To formulate the main result Theorem, we need some concepts related to polytopes and fans (see \cite{Cx}, \cite{Zg}).
First of all, remark that the   hyperplanes mentioned above divide  $\mathbb{R}^n$ into cones that form a fan (see  \cite{Zg}, 7.1).
By definition, a fan is a conic polyhedron $\Sigma\subset\mathbb{R}^n.$ This means that $\Sigma$ is a family of strongly convex rational cones $\{\sigma\}$ with the properties:

$ (i) \;\; \text{if} \;\;\sigma \in \Sigma \;\; \text{then all faces of} \;\; \sigma \;\; \text{belong to }  \Sigma; $

$ (ii) \;\; \text{if} \;\;\sigma,\tau \in \Sigma \;\; \text{then } \;\; \sigma \cap \tau \in \Sigma. $\\
Denote by  $\mu_{1},...,\mu_{d} $ the one-dimensional generators of this fan. They define the polar (with respect to  the fan) polytope (see ~\cite{Zg}, Corollary 7.18) 
$$P=\{\alpha\in\mathbb{R}^n: (\mu_\nu,\alpha)\geq {g}(\mu_\nu) \}, \,\, \nu=1,...,d.$$
Note that the polytope depends on the choice of $\varepsilon_p$ in $g$, and in some cases the interior of $P$ may be empty. We are interested in the case when it is not empty.\\

\textbf{Theorem. } \textit{Assume that the interpolating function } $\varphi(\zeta)$ \textit{is holomorphic in   } $\mathbb{R}_{\geq0}^n + i\mathbb{R}^n$ \textit{and satisfies the estimates}\\

$ 1) \;\log|\varphi(\xi+i\eta)|\leq\tilde g(\eta) \;\;
\textit{on} \;\; \textit{imaginary subspace} \;\; i\mathbb{R}^n, $\\

$ 2)\; \log{|\varphi(\xi+i\eta)|}\leq \overset{n}{\underset{j=1}\sum}(\pi-\delta)|\eta_j|+b\xi_j+C \;\; \textit{on}  \; \; \mathbb{R}_{\geq0}^n + i\mathbb{R}^n,   $\\
\textit{for some constants $b$,  $C,$ and $\delta>0.$}\\
\textit{Then the sum of the series $(\ref{1.0})$ extends analytically into the sectorial domain} $Arg^{-1}(P^o),$ \textit{where} $P^o$ \textit{is the interior of the }  \textit{polyhedron} $P$. \\ 

Note that this theorem essentially improves our previous result, Theorem 1 of \cite{Mk2}. Recall that there we required  the interpolating function to be holomorphic on  $(-\delta, \infty)^n + i\mathbb{R}^n$ for some $\delta>0$,   condition 1) to hold on $   (-{\delta},0]^n+i\mathbb{R}^n, $ and condition 2) to hold on $ (-\delta, \infty)^n\setminus {(-{\delta},0]^n} + i\mathbb{R}^n.$
\\

In order to prove this theorem, we  need to reformulate and strengthen the Le Roy-Lindel\"of theorem for power series in one variable. 

Consider a  power series 
\begin{equation}\label{10.0} 
{\sum_{k=0}^\infty f_{k} {z}^{k}    }   
\end{equation}
with radius of convergence  $R$. According to the Cauchy-Hadamard theorem this means that 
\begin{equation*}\label{10.1}
{\varlimsup_{n \to \infty}{ ^{n}\sqrt{|f_{n}|}=\frac{1}{R}.}}
\end{equation*}

\textbf{Proposition } \textit{Assume that the  function} $\varphi(\zeta)$ \textit{is holomorphic in   } $\mathbb{R}_{\geq0} + i\mathbb{R}$ \textit{and interpolates the coefficients of power series $(\ref{10.0})$. If $\varphi(\zeta)$ satisfies the condition}

\begin{equation}\label{T2.1}
 \max\left\{h_{\varphi}(\frac{\pi}{2}), h_{\varphi}(-\frac{\pi}{2})\right\} = \sigma < \pi, 
\end{equation}
\textit{then the sum of the series (\ref{10.0})   extends analitically to the sector} $\mathbb{C}\setminus \Delta_\sigma$. \\

Proof of the Proposition.\\

Recall that $h_{\varphi}(\theta)$ is  \textit{trigonometric convex}\cite{Bo}
 for $\alpha_1<\theta<\alpha_2$ and  $\alpha_2-\alpha_1<\pi,$  the following inequality holds
\begin{equation*}\label{trig_convex}
	h_{\varphi}(\theta)\sin{(\alpha_2-\alpha_1)}\leq h_{\varphi}(\alpha_1)\sin{(\alpha_2-\theta)}+h_{\varphi}(\alpha_2)\sin{(\theta-\alpha_1)}.
\end{equation*}

Appling this property to intervals $[0, \frac{\pi}{2}]$,  $[-\frac{\pi}{2}, 0]$ and using $(\ref{T2.1})$ we get 

\begin{equation}\label{Ind1}
h_{\varphi}(\theta) \leq h_{\varphi}(0)\cos(\theta)+\sigma|\sin(\theta)| \;\;\; \textit{for} \;\;\; |\theta|\leq \frac{\pi}{2}.
\end{equation}

According to  Bernstein`s theorem (\cite{Levin}, p.100)

\begin{equation}\label{}
\limsup_{n \to \infty}  \frac{\ln\left|\varphi(n)\right|}{n}=
\limsup_{r \to \infty}  \frac{\ln\left|\varphi(r)\right|}{r},
\end{equation}
hence 
\begin{equation}\label{Rad}
R= e^{-h_{\varphi}(0)}.
\end{equation}

The definition $(\ref{indicator})$ yields
\begin{equation*}{\label{1.11}}
{|\varphi(re^{i\theta})| \leq e^{h_\varphi (\theta)r+o(r)} \;\; \textit{for} \;\; |\theta|\leq \frac{\pi}{2} \;\; \textit{as} \;\; r\to \infty.
}
\end{equation*}
Therefore using $(\ref{Ind1})$ and $(\ref{Rad})$ we get 
%\begin{equation} {\label{11.}}
%{|\varphi(\zeta)|<ce^{h_{\varphi}(0) r \cos(\theta)+\sigma r |\sin(\theta)|+o(r)}.
%}
%\end{equation}

\begin{equation} {\label{11.}}
{|\varphi(\zeta)|<ce^{\ln R^{-1} r \cos(\theta)+\sigma r |\sin(\theta)|+o(r)}\;\; \textit{for} \;\; |\theta|\leq \frac{\pi}{2}.
}
\end{equation}

Introduce an auxiliary function
$$ g(\zeta,z):= \frac{z^{\zeta} }{e^{2\pi i \zeta}-1 }.$$
Let $D_{\frac{1}{8}}(a):=\{\zeta\in \mathbb{C}: |\zeta-a|<\frac{1}{8}\}$  be a disk of radius $\frac{1}{8}$  centered at  $a\in \mathbb{C}$.  Denote $$D^*:=\Delta_{\frac{\pi}{2}}\setminus(\cup_{m\in \mathbb{Z}}D_{\frac{1}{8}}(m)) \;\; \text{and } \;\; K=\bar D_{Re^{-\delta}}\setminus(\Delta_{\sigma+\delta}^o\cup D_{\delta}).$$
Note that there exists a constant $C>0$ such that 
\begin{equation*} \label{d1.3}
{ |e^{2\pi i \zeta}-1|>\frac{e^{\pi(|\mathop{\text{Im}} \zeta|-\mathop{\text{Im}} \zeta)}}{C} \;\;\; \text{ for} \;\;\; \zeta \in D^*.
}
\end{equation*}
Using this  inequality we get the estimation
\begin{equation}\label{d1.4}
{ |g(\zeta,z)|<Ce^{r\cos\theta \ln|z|- (\pi-|\pi-\arg z|) r|\sin\theta|}
}
\end{equation}
for $\zeta \in D^*$  and
$ z \in \mathbb{C} \setminus \mathbb R_+.$\\
%Let
%\begin{equation*}
%{K=\bar D_{R-\delta}\setminus(\Delta_{\sigma+\delta}^o\cup D_{\delta}).
%}
%\end{equation*}
Thus, for $\zeta \in D^*$  and
$ z \in \mathbb{C} \setminus \mathbb R_+$
\begin{equation} {\label{11.}}
{|\varphi(\zeta)||g(\zeta,z)|<ce^{r\cos \theta( \ln|z|-\ln R)-(\pi-\sigma-|\pi-\arg z|)r|\sin\theta|+o(|r|)}.
}
\end{equation}

Note that
$ \pi-\sigma-|\pi-\arg z|\geq\delta \;\; \text{ for}\;\; z\in \mathbb{C}\setminus \Delta_{\sigma+\delta}$
 and 
$$\ln|z|-\ln R\leq -\delta\;\; \text{ for}\;\;z\in K:=\bar D_{Re^{-\delta}}\setminus(\Delta_{\sigma+\delta}^o\cup D_{\delta}).$$
Hence, for  $z\in {K} \; $   and  $ \; \zeta \in D^*$
we get
\begin{equation} {\label{11.}}
{|\varphi(\zeta)||g(\zeta,z)|<ce^{-\delta r(\cos \theta+|\sin\theta|)+o(|r|)}.
}
\end{equation}

For any  $m\in\mathbb{N}$ we consider the integral
\begin{equation} {\label{2.6}}
{I_m=\int\limits_{\partial G_{m}}\varphi(\zeta)g(\zeta,z)d\zeta},
\end{equation}
	where $G_{m}$ is a planar domain 
bounded by a deformed segment
	$$\Gamma_{m}^{1}=\left[-i(m+\frac{1}{2}),-{i\frac{1}{2}}\right]\cup\left\{\frac{1}{2} e^{i\theta},|\theta|\leq\frac{\pi}{2}\right\} \cup\left[{i\frac{1}{2}},i(m+\frac{1}{2})\right],$$
		and a half circle 
		$$\Gamma_{m}^{2}=\left\{(m+\frac{1}{2})e^{i\theta},|\theta|\leq\frac{\pi}{2}\right\}.$$
Represent the integral $I_m$ as the sum of integrals
 $I_m^1$ and  $ I_m^2$  over
 $\Gamma_{m}^1$ and $ \Gamma_{m}^2, $ respectively.

For $ \zeta\in  D^* $ and $z\in K$ we obtain the estimate 
\begin{equation*}
I_m^2= \int\limits_{\Gamma_m^2}|\varphi(\zeta)g(\zeta,z)||d\zeta|\leq c(m+\frac{1}{2})e^{-\delta (m+\frac{1}{2})(\cos\theta+|\sin\theta|)}\leq \acute{c}(m+\frac{1}{2})e^{-\delta (m+\frac{1}{2})}.
\end{equation*}
Thus, for $z\in K$ integral   $I_m^2 $ tends to zero as $m\to \infty$ and
\begin{equation*}{\label{1.20}}
\lim_{m\to\infty}I_m =\lim_{m\to\infty}\int\limits_{\partial G_m}\psi(\zeta)g(\zeta,z)d\zeta=\lim_{m\to\infty}\int\limits_{\Gamma_m^1}\psi(\zeta)g(\zeta,z)d\zeta=\lim_{m\to\infty}I_m^1.
\end{equation*}

In the domain  $G_m$ the integrand has simple poles at real integer points.\\
The residue theorem yields 
\begin{equation*}{\label{1.21}}
\int\limits_{\partial G_m}\varphi(\zeta)g(\zeta,z)d\zeta=\sum_{n=1}^{m}\varphi(n)z^n.
\end{equation*}

Consider the integral
\begin{equation*}{\label{1.22}}
I=\int\limits_{\Gamma^{1}}\varphi(\zeta)g(\zeta,z)d\zeta.
\end{equation*}

$$\text{Denote by } \;\; \Gamma^{1}=\left[-i\infty,-{i\frac{1}{2}}\right]\cup\left\{\frac{1}{2} e^{i\theta},|\theta|\leq\frac{\pi}{2}\right\} \cup\left[{i\frac{1}{2}},i\infty\right]  \;\;\; \text{the deformed imaginary axis}.$$
For $ \zeta=i\eta $ and $z\in \mathbb{C}\setminus (\Delta_{\sigma+\delta}\cup D_{\delta})$ using $(\ref{11.})$ we find

\begin{equation*} {\label{1.23}}
{|\varphi(\zeta)||g(\zeta,z)|<ce^{-\delta|\eta|+o(|\zeta|)}.
}
\end{equation*}

It follows from this inequality that the integral  $I$ converges absolutely and uniformly on any compact subset $  K \subset  \mathbb{C}\setminus (\Delta_{\sigma+\delta}\cup D_{\delta}) , $ and defines a holomorphic function on the set of interior points of $
K^o$.
For  $K$
\begin{equation*}{\label{1.24}}
\int\limits_{\Gamma_m^1}\varphi(\zeta)g(\zeta,z)d\zeta\to\int\limits_{\Gamma^1}\varphi(\zeta)g(\zeta,z)d\zeta, \;\;  \text{as} \;\;  m\to \infty.
\end{equation*}

Since $I_m \to  I \;\; \text{as} \;\; m\to \infty,$
 $ I(z)=f(z)$  for $z \in D_R\cap K^o$. This means that  $f(z)$  extends analytically to $K^o$. Since  $K$ is an arbitrary compact set in $\mathbb{C}\setminus \Delta_{\sigma+\delta}$ for any small  $\delta$,   the function   $f(z)$ extends to the open sector $\mathbb{C}\setminus \Delta_{\sigma}$. \\ 

	Examples.
	
	Consider the power series 
	\begin{equation*}\label{1p.1}
{\sum_{k=0}e^kz^k, \;\;\; R=e^{-1}
}
\end{equation*}
  representing the rational function
$$f(z)=\frac{1}{1-ez}.$$ 
Obviously, the function $\varphi(\zeta)=e^z$ interpolates the coefficients of the series and holomorphic in $\mathbb{R}_{\geq0} + i\mathbb{R}.$
The exponential type  $\sigma$ of  the function $\varphi(\zeta)$ in $\mathbb{R}_{\geq0} + i\mathbb{R}$ equals 1 and according to the Le Roy-Lindel\"of theorem the series extends to the sector   $\mathbb{C}\setminus \Delta_{1}$.\\
The indicator function is  
$$h_{\varphi}(\theta)=\cos\theta,$$
hence  $$\max\left\{h_{\varphi}(\frac{\pi}{2}), h_{\varphi}(-\frac{\pi}{2})\right\} = 0 < \pi.$$
Therefore, according to the   proposition the series extends to the following sector $\mathbb{C}\setminus \Delta_{0}$. Note that the sector $\mathbb{C}\setminus \Delta_{0}$ contains sector $\mathbb{C}\setminus \Delta_{1}$ and it is the maximal sector where  this series extends.
Also note that here we can not apply  Arakelian's theorem $\cite{Ar84}$ because the radius of convergence may be not equal to 1.
The proof of theorem. The reasoning is similar to that
 in the article \cite{Mk2}, but there is a difference and for convenience we  give the complete proof.

Consider  an auxiliary function
$$ h(\zeta,z)=\prod_{j=1}^n \frac{z_j^{\zeta_j} }{(e^{2\pi i \zeta_j}-1) }.$$
The function $h(\zeta,z)$  is meromorphic in   $\zeta \in \mathbb{C}^n$ and holomorphic in $z \in (\mathbb{C} \setminus \mathbb R_+)^n.$\\
Let us introduce  the following product
\begin{equation*}\label{d1.14}
{ \prod_{j=1}^n \frac{z^{\zeta_j} }{(e^{2\pi i \zeta_j}-1)}\varphi(\zeta),
}
\end{equation*}
it is  meromorphic in   
$\mathbb{R}^n_{\geq 0}~+~i\mathbb{R}^n$ with poles only on divisors
\begin{equation*}\label{d1.15}
{ Q_1=\{(\zeta_1, ... ,\zeta_n) :  e^{2\pi i \zeta_1}-1=0\}=\mathbb{Z}\times \mathbb{C}^{n-1} ,
}
\end{equation*}
\begin{equation*}\label{d1.16}
{ .\; . \; . \; . \;  . \; . \; . \; . \; . \;  . \; . \; . \;  . \;  . \; . \;  . \; . \;  . \;  . \; . \;  . \;  . \;  . \;  .
}
\end{equation*}
\begin{equation*}\label{d1.17}
{ Q_n=\{(\zeta_1, ... ,\zeta_n) : e^{2\pi i \zeta_n}-1=0\}=\mathbb{C}^{n-1} \times  \mathbb{Z}.
}
\end{equation*}
For each vector $m=(m_1,...,m_n)\in\mathbb{N}^n$ we define a domain $G_m=G_{m_1}\times ... \times  G_{m_n}$  where $G_{m_j}$ is a planar domain in $\mathbb{C}_{\zeta_j}$  bounded by a deformed segment
$$\Gamma_{m_j}^{1}=[-i(m_j+\frac{1}{2}),-{i\frac{\delta}{2}}]\cup\{\frac{\delta}{2} e^{i\theta_j},|\theta_j|\leq\frac{\pi}{2}\} \cup[{i\frac{\delta}{2}},i(m_j+\frac{1}{2})]$$
		and a half circle 
		$$\Gamma_{m_j}^{2}=\{(m_j+\frac{1}{2})e^{i\theta_j},|\theta_j|\leq\frac{\pi}{2}\}.$$
		$G_m$ is compatible with the family of divisors $Q_1, ... ,Q_n$ in the sense that for each $j\in\{1,...,n\}$ the edge
$$\bar{G}_{m_1}\times ... \times\bar{G}_{m_{j-1}} \times \partial\bar{G}_{m_j}\times \bar{G}_{m_{j+1}}\times ... \times\bar{G}_{m_{n}}$$ does not intersect the divisors 
$$Q_j=\mathbb {C}\times ...\times \mathbb {C}\times\underset{j}{\mathbb {Z}}\times \mathbb {C} \times ... \times\mathbb {C}.$$
With this compatibility property the residue theorem  (see \cite{Ts},\cite{TJ}) states that
%we get the following representation for partial sums of the series (1):
\begin{equation} {\label{13}}
{\frac{1}{2\pi i}\int\limits_{\partial_0 G_{m}}\varphi(\zeta)h(\zeta,z)d\zeta=\sum_{k_1=1}^{m_1}... \sum_{k_n=1}^{m_n}\varphi(k_1,...,k_n)z_1^{k_1}...z_n^{k_n},
}
\end{equation}
where $\partial_0 G_{m}=\partial G_{m_1}\times ... \times \partial G_{m_n}$ is the distinguished boundary of the domain $G_{m}$.

By condition 2) of the Theorem  we get 
\begin{equation}{\label{d1.2}}
{|\varphi(\zeta)|\leq A e ^{(\pi-\delta)\sum |\eta_j|+b\sum \xi_j}  \;\;\; \text{for} \;\;\; \zeta \in  \mathbb{R}_{\geq0}^n + i\mathbb{R}^n.}
\end{equation}
%for $\zeta \in \Delta_\delta:= (-\delta, \infty)^n\setminus {(-{\delta},0]^n} + i\mathbb{R}^n.$

We get the estimation
\begin{equation}\label{d1.4}
{ |h(\zeta,z)|<C_{\delta}e^{\left\langle \xi,\ln|z|\right\rangle-\left\langle (\pi-|\pi-\arg z|),|\eta|\right\rangle}
}
\end{equation}
for
$ z \in (\mathbb{C} \setminus \mathbb R_+)^n$ and  $\zeta \in (\mathbb{C}\setminus D)^n$ where $D:=(\cup_{m\in \mathbb{Z}}D_{\frac{\delta}{8}}(m))$. \\
Combining (\ref{d1.2}) and (\ref{d1.4}) for $$(\zeta,z)\in ((\mathbb{R}_{\geq0} + i\mathbb{R})\cap D)^n \times (\mathbb{C} \setminus \mathbb R_+)^n$$ we obtain
\begin{equation*}\label{d1.5}
{|\varphi(\zeta)||h(\zeta,z)|\leq c e^{\sum (\xi_j\ln|z_j|+\xi_jb)-\sum (\delta-|\pi-\arg z_j|)|\eta_j|}
}.
\end{equation*}

%\begin{minipage}{0.6\linewidth}
%	Denoting $d(z_j)=\frac{\delta}{2}-|\pi-\arg z_j|$, we rewrite this inequality in the form
%	\begin{equation*}\label{d1.7}
%	{|\varphi(\zeta)||h(\zeta,z)|\leq c e^{\sum (\xi_j\ln|z_j|+|\xi_j|b)-\sum d(z_j)|\eta_j|}.
%	}
%	\end{equation*}
	Note  that $\delta-|\pi-\arg z_j|\geq\frac{\delta}{2} $ when $z_j$ runs over the compact 
	\begin{equation*}\label{d1.8}
	{K=(\bar D_{e^{-(b+\delta)}}\setminus D_{e^{-2(b+\delta)}}) \setminus\Delta_{\pi-\frac{\delta}{2}}^o.
	}
	\end{equation*}
	Finally, for  $z\in {K}^n$ \; and  \; $\zeta \in((\mathbb{R}_{\geq0} + i\mathbb{R})\cap D)^n$
	we arrive at the estimation
\begin{equation}\label{d1.10}
{ |h(\zeta,z)||\varphi(\zeta)|<c e^{-\delta\sum \xi_j-\frac{\delta}{2} \sum |\eta_j|}.
}
\end{equation}

Denote the partial sum in $(\ref{13})$ by  $I_m$. Due to the construction $G_{m_j}$, we get  by  $(\ref{13})$ that  $I_m$ can be represented as a sum of $2^n$ integrals over  pieces
$$ \Gamma_{m_1}^{i_1}\times ... \times \Gamma_{m_n}^{i_n}, \;\;\; i_k\in{1,2}.$$

For each such piece we split the integration variables $\zeta_1,...,\zeta_n$  into two groups: $B_1$ and $ B_2,$ where $j \in B_l$ if $\zeta_j \in \Gamma_{m_j}^{l}.$
When $z_j\in K$ and $j\in {B_1}$, i.e. 
$$\zeta_j\in\Gamma_{m_j}^{1}=[-i(m_j+\frac{1}{2}),-{i\frac{\delta}{2}}]\cup\{\frac{\delta}{2} e^{i\theta_j},\frac{\pi}{2}\leq\theta_j\leq\frac{3\pi}{2}\} \cup[{i\frac{\delta}{2}},i(m_j+\frac{1}{2})],$$ 
we get 
$$e^{({\xi}_j\ln|z_j|+|{\xi}_j|b)-\frac{\delta}{2}|{\eta}_j|}\leq e^{c_j-\frac{\delta}{2}|{\eta}_j|} = c^{'}_j 
e^{-\frac{\delta}{2}|{\eta}_j|}.$$

When $z_j\in K$ and $j\in {B_2}$, i.e. $$\zeta_j\in\Gamma_{m_j}^{2}=\{(m_j+\frac{1}{2})e^{i\theta_j},|\theta_j|\leq\frac{\pi}{2}\},$$ 
we get
$$e^{(\xi_j\ln|z_j|+|\xi_j|b)-\frac{\delta}{2}  |\eta_j|}\leq e^{(\xi_j(\ln|z_j|+b)-\frac{\delta}{2}|\eta_j|} \leq e^{-\frac{\delta}{2}(\xi_j+|\eta_j|)}\leq c_je^{-\frac{\delta}{2}m_j}.$$
Then using ($\ref{d1.10}$) for $z\in K^n$  we obtain the following estimate

\begin{equation*} \label{d1.12}
{|h(\zeta,z)||\varphi(\zeta)|<c \prod_{j\in {B_1}} e^{-\delta|\eta_j|} \prod_{j\in B_2} e^{-\delta m_j}.}
\end{equation*}

Thus, if  $B_2\neq \emptyset$,  the integral over the corresponding contour vanishes as $m_j\to \infty, \;\; j=1,...,n.$ 
Therefore, for  $z\in K^n $  we get
\begin{equation*} \label{d1.13}
{I_m=\frac{1}{(2\pi i)^n}\int\limits_{\partial G_{m}}\varphi(\zeta)h(\zeta,z)d\zeta \to \frac{1}{(2\pi i)^n}\int\limits_{ \Gamma_m}\varphi(\zeta)h(\zeta,z)d\zeta \;\; \text{as} \;\; m_j \to \infty }
\end{equation*}
where $\Gamma_m=\Gamma_{m_1}^{1}\times ... \times \Gamma_{m_n}^{1}.$

Consider the integral

\begin{equation*}\label{d1.23}
{ I=\int\limits_{\Gamma} h(\zeta,z)\varphi(\zeta)d\zeta,
}
\end{equation*}
where $$\Gamma=\left((-i\infty, -i\delta)\cup [\delta e^{i\theta}, |\theta|\leq\frac{\pi}{2}]\cup (i\delta, i\infty)\right)^n$$ is the deformed imaginary subspace $i\mathbb{R}^n.$

We obtain that the absolute value of the integrand $|\varphi(\zeta)||h(\zeta,z)|$ is estimated by
\begin{equation*}\label{d1.24}
|h(\zeta,z)||\varphi(\zeta)|<c \exp \{\sum_{p=1}^{q}\varepsilon_p|a_p^1\eta_1+...+a_p^n\eta_n|+\pi\sum_{j=1}^n(\eta_j-|\eta_j|)-\sum_{j=1}^n\eta_j\arg z_j\}.   
\end{equation*}

For the convergence of the integral (see p.164 \cite{ST}) it is necessary that the following inequality holds 
\begin{equation*}\label{d1.25}
{\sum_{j=1}^n\eta_j\arg z_j}>{\sum_{p=1}^{q}\varepsilon_p|a_p^1\eta_1+...+a_p^n\eta_n|+\pi\sum_{j=1}^n(\eta_j-|\eta_j|),
}
\end{equation*}

which is the same as 
\begin{equation}\label{d1.27}
{(\eta,\arg z)> g(\eta).
}
\end{equation}
The inequality $(\ref{d1.27})$ for $\alpha=\arg z$ holds for interior points of the polyhedron  
$$P=\{\alpha\in\mathbb{R}^n: (\mu_\nu,\alpha)\geq g(\mu_\nu) \}$$
where  $\pm\mu_{1},...,\pm\mu_{d} $ are vectors of  $\mathbb{R}^n$ generating the fan defined by the hyperplanes $a_p^1\eta_1+...+a_p^n\eta_n~=~0,$ \,\, $p~=~1,...,q$  \; and \; $\eta_k=0, \,\, k=1,...,n.$ \\
Thus, the integral $\mathbb I$ converges for  $z\in Arg^{-1}(P^o)$.

Since $I_m \to  I\;\; \text{as} \;\; m_j\to \infty,\;\; j=1,...,n,$
we get 
$$ f(z)=I(z)+\sum_{l=1}^nf_l(z)$$
  for  $z \in  (K^o)^n$, where 
	$$f_l(z)=\sum_{k\in\mathbb N^n, k_l=0} f_{k} {z}^{k}.$$
	Now we show that the functions $f_l(z)$ are analytic for $z\in Arg^{-1}(P^o)$.
	We shall show this in the case when $n=2$  the general case is proved by induction.
So, for $n=2$ we have $$ f(z_1,z_2)=I(z_1, z_2)+\sum_{k_1=0}^\infty f_{k_10} {z_1}^{k_1} + \sum_{k_2=0}^\infty f_{0k_2} {z_2}^{k_2}        $$  for  $z \in  (K^o)^n$.
%Therefore the sum of the series (\ref{1.0})  extends analytically into the sectorial domain $Arg^{-1}(P^o),$  which  was to be proved.

Consider the series 
\begin{equation}\label{s2.1}
{\sum_{k_1=0}^\infty f_{k_10} {z_1}^{k_1}
}
\end{equation} 
The function $\varphi_1(\zeta_1):=\varphi(\zeta_1, 0)$ interpolates the coefficients of the power series $(\ref{s2.1})$
Denote 
$$\sigma_1:=\sum_{p=1}^{q}\varepsilon_p|a_p^1|.$$ 
According to the construction of set $P$ we obtain that 
$\sigma_1\leq \alpha_1 \;\;  \text{ for any} \;\; \{\alpha_1\in\mathbb{R} : (\alpha_1, \alpha_2)\in P\}.$ 
Note that $\max\{h_{\varphi_1}(\frac{\pi}{2}), h_{\varphi_1}(-\frac{\pi}{2})\}\leq\sigma_1<\pi$ 
Hence, according to the Proposition the series $(\ref{s2.1})$ extends into the sector
$$ \mathbb{C}\setminus \Delta_{\sigma_1}$$
 in terms of argument $\sigma_1<\arg z_1< 2\pi-\sigma_1.$
As a function of the variables $(z_1, z_2)$ it extends continue into the argument domain  $\sigma_1<\arg z_1< 2\pi-\sigma_1, 0<\arg z_2< 2\pi$ which includes  $Arg^{-1}(P^o),$
Similarly, the series $$\sum_{k_2=0}^\infty f_{0k_2} {z_2}^{k_2}        $$ extends into a  domain which includes  $Arg^{-1}(P^o).$

Therefore the sum of the series (\ref{1.0})
$$ f(z_1,z_2)=I(z_1, z_2)+\sum_{k_1=0}^\infty f_{k_10} {z_1}^{k_1} + \sum_{k_2=0}^\infty f_{0k_2} {z_2}^{k_2}$$  extends analytically into the sectorial domain $Arg^{-1}(P^o),$  since each function of the right-hand side of the last equality continues into the sectorial domain $Arg^{-1}(P^o)$.

\vbpb
\end{document}